\documentclass[reqno]{amsart}
\usepackage{url}
\newtheorem{theorem}{Theorem}

\theoremstyle{remark}

\allowdisplaybreaks
\begin{document}
\title[A necessary condition for the boundedness]
{A necessary condition for the boundedness of the maximal operator on 
$L^{p(\cdot)}$ over reverse doubling spaces of homogeneous type}
%%%----------------------------------------------------------------------------
\author[O.~Karlovych]{Oleksiy Karlovych}  
\address{
Centro de Matem\'atica e Aplica\c{c}\~oes,
Departamento de Matem\'a\-tica,
Faculdade de Ci\^en\-cias e Tecnologia,
Universidade Nova de Lisboa,
Quinta da Torre,
2829--516 Caparica, Portugal}
\email{oyk@fct.unl.pt}
%%%----------------------------------------------------------------------------
\author[A.~Shalukhina]{Alina Shalukhina}  
\address{
Centro de Matem\'atica e Aplica\c{c}\~oes,
Departamento de Matem\'a\-tica,
Faculdade de Ci\^en\-cias e Tecnologia,
Universidade Nova de Lisboa,
Quinta da Torre,
2829--516 Caparica, Portugal}
\email{a.shaukhina@campus.fct.unl.pt}
%%%----------------------------------------------------------------------------
\thanks{
This work is funded by national funds through the FCT - Funda\c{c}\~ao para a 
Ci\^encia e a Tecnologia, I.P., 
under the scope of the projects UIDB/00297/2020 
(\url{https://doi.org/10.54499/UIDB/00297/2020})
and UIDP/00297/2020 
(\url{https://doi.org/10.54499/UIDP/00297/2020})
(Center for Mathematics and Applications).
The second author is funded by national funds through the FCT – 
Funda\c{c}\~ao para a Ci\~encia e a Tecnologia, I.P., 
under the scope of the PhD scholarship UI/BD/154284/2022.
}
%%%----------------------------------------------------------------------------
\begin{abstract}
Let $(X,d,\mu)$ be a space of homogeneous type and $p(\cdot):X\to[1,\infty]$
be a variable exponent. We show that if the measure $\mu$ is Borel-semiregular
and reverse doubling, then the condition 
$\operatornamewithlimits{ess\,inf}\limits_{x\in X}p(x)>1$ 
is necessary for the boundedness of the Hardy-Littlewood maximal operator $M$ 
on the variable Lebesgue space $L^{p(\cdot)}(X,d,\mu)$. 
\end{abstract}
%%%----------------------------------------------------------------------------
\maketitle 
%%%----------------------------------------------------------------------------
\section{Introduction and the main result}
Let $(X,d,\mu)$ be a space of homogeneous type 
(see Section~\ref{sec:preliminaries}). For $x\in X$ and $r>0$, 
consider the ball $B(x,r):=\{y\in X:d(x,y)<r\}$ centered at $x$ of radius $r$.
By definition, the Borel measure $\mu$ has the doubling property, that is, 
for every ball $B(x,r)$ with respect to the quasi-metric $d$, the relation
\[
\mu(B(x,r))\le A\mu(B(x,r/2))
\]
holds true with an absolute constant $A>1$. We will assume that 
$0<\mu(B)<\infty$ for every ball $B$. If, additionally, the 
reverse inequality
\[
\mu(B(x,r/2))\le \delta\mu(B(x,r))
\] 
is valid for some absolute constant $0<\delta<1$, it is said that 
the measure $\mu$ is also reverse doubling (see, e.g., \cite[Section~2.6]{DFGL08}).

Given a complex-valued function $f\in L^1_{\rm loc}(X,d,\mu)$, we define its 
Hardy-Littlewood maximal function $Mf$ by 
%%%
\begin{equation}\label{eq:Hardy-Littlewood}
Mf(x):=\sup_{B\ni x}\frac{1}{\mu(B)}\int_B |f(y)|\,d\mu(y),
\quad x\in X,
\end{equation}
%%%
where the supremum is taken over all balls $B\subset X$ containing $x\in X$.
The Hardy-Little\-wood maximal operator $M$ is a sublinear operator acting by
the rule $f\mapsto Mf$.

Let $L^0(X,d,\mu)$ denote the set of all complex-valued measurable functions
on $X$ and let $\mathcal{P}(X)$ denote the set of all measurable functions
$p(\cdot):X\to[1,\infty]$. The functions in $\mathcal{P}(X)$ are called
variable exponents. Let $X_\infty:=\{x\in X:p(x)=\infty\}$.
For a function $f\in L^0(X,d,\mu)$ and $p(\cdot)\in\mathcal{P}(X)$, consider 
the functional, which is called modular, given by
\[
\varrho_{p(\cdot)}(f)
:=
\int_{X\setminus X_\infty} |f(x)|^{p(x)}\,d\mu(x)
+
\operatornamewithlimits{ess\,sup}_{x\in X_\infty}|f(x)|.
\]

By definition, the variable Lebesgue space $L^{p(\cdot)}(X,d,\mu)$ consists 
of all functions $f\in L^0(X,d,\mu)$ such that 
$\varrho_{p(\cdot)}(f/\lambda)<\infty$
for some $\lambda>0$ depending on $f$. It is a Banach space with respect
to the Luxemburg-Nakano norm given by
\[
\|f\|_{p(\cdot)}:=\inf\{\lambda>0\ :\ \varrho_{p(\cdot)}(f/\lambda)\le 1\}.
\]
If $p(\cdot)\in\mathcal{P}(X)$ is constant, then $L^{p(\cdot)}(X,d,\mu)$ is 
nothing but the standard Lebesgue space $L^p(X,d,\mu)$. Variable Lebesgue spaces 
are often called Naka\-no spaces. We refer to Maligranda's paper \cite{M11} 
for the role of Hidegoro Nakano in the study of variable Lebesgue spaces
and to the monographs \cite{CF13,DHHR11} for the basic properties of these
spaces. 

Cruz-Uribe, Fiorenza and Neugebauer proved in \cite[Theorem~1.7]{CFN03} 
that if $\Omega\subset\mathbb{R}^d$ is an open set and 
$p(\cdot):\Omega\to[1,\infty)$ is upper semi-continuous, then the boundedness of
$M$ on $L^{p(\cdot)}(\Omega)$ implies that $\inf\limits_{x\in\Omega}p(x)>1$
(it is supposed that $\Omega$ is equipped with the Lebesgue measure
and the usual Euclidean distance).
The upper semi-continuity assumption was removed by
Deining et al. \cite[Theorem~1.6]{DHHMS09}. Another proof of this fact was
given by Izuki, Nakai and Sawano in \cite[Proposition~3.3]{INS13}
and \cite[Proposition~21.2]{ISN14}. Proofs of the fact that the 
boundedness of $M$ on $L^{p(\cdot)}(\mathbb{R}^d)$ implies that 
\[
p_-(\mathbb{R}^d):=
\operatornamewithlimits{ess\,inf}\limits_{x\in\mathbb{R}^d}p(x)>1
\]
are given in \cite[Theorem~4.7.1]{DHHR11} and \cite[Theorem~3.19]{CF13}.
Note also that very recently Roberts \cite[Theorem~3.3]{R23} extended the 
above result to the setting of the fractional maximal operator
\[
M_\alpha f(x):=\sup_{Q\ni x}|Q|^{\alpha/d-1}\int_Q|f(y)|dy,
\]
where $0\le \alpha<d$ and the supermum is taken over all cubes
$Q\subset\mathbb{R}^d$ containing $x$, and proved that 
if $M_\alpha$ is bounded from
$L^{p(\cdot)}(\mathbb{R}^d)$ to $L^{q(\cdot)}(\mathbb{R}^d)$ with
$p(\cdot)\in\mathcal{P}(\mathbb{R}^d)$ and $q(\cdot)$ defined by
$1/p(\cdot)-1/q(\cdot)=\alpha/d$, then $p_-(\mathbb{R}^d)>1$.

In August of 2019, during the ISAAC Congress held in Aveiro, Portugal,
Stefan Samko asked the first author, under which conditions on 
a quasi-metric space $(X,d)$ and a measure $\mu$ on $X$, the 
boundedness of the Hardy-Littlewood maximal operator $M$ on 
$L^{p(\cdot)}(X,d,\mu)$ implies that
\[
p_-(X):=\operatornamewithlimits{ess\,inf}_{x\in X} p(x)>1.
\] 
Surprisingly enough, we were not able to find any result on necessary
conditions for the boundedness of $M$ on variable Lebesgue spaces 
beyond the Euclidean setting. Moreover, we are not aware of a proof
that $M$ is unbounded on $L^1(X,d,\mu)$ in the setting of spaces of
homogeneous type.

The aim of this paper is to address this open problem. Our main result is 
the following.
%%%----------------------------------------------------------------------------
\begin{theorem}\label{th:main}
Suppose $(X,d,\mu)$ is a space of homogeneous type which has the property that 
the measure $\mu$ is Borel-semiregular and reverse doubling. Given an exponent 
function $p(\cdot)\in\mathcal{P}(X)$, if the Hardy-Littlewood maximal operator
$M$ is bounded on the variable Lebesgue space $L^{p(\cdot)}(X,d,\mu)$, then
$p_-(X)>1$.
\end{theorem}
%%%----------------------------------------------------------------------------
Although the assumption that the measure $\mu$ is reverse doubling is
essential in our proof, we believe that Theorem~\ref{th:main} should be 
true without it.

The paper is organized as follows. In Section~\ref{sec:preliminaries},
we provide necessary background on spaces of homogeneous type and the
Lebesgue differentiation theorem in this setting. Section~\ref{sec:proof}
contains the proof of Theorem~\ref{th:main}.
We conclude this paper observing in Section~\ref{sec:final-remark}
that the hypothesis of Borel-semiregularity
of $\mu$ can be dropped if one requires that the variable exponent
$p(\cdot)$ is upper semi-continuous.
%%%----------------------------------------------------------------------------
\section{Preliminaries on spaces of homogeneous type}\label{sec:preliminaries}
Following \cite[Section~2.1]{AM15}, given a nonempty set $X$, call a function
$\varrho:X\times X\to[0,\infty)$ a quasi-distance (or a quasi-metric)
provided there exist constants $C_0,C_1\in(0,\infty)$ such that for 
all $x,y,z\in X$ the following axioms hold:
\begin{enumerate}
\item[(a)] $\varrho(x,y)=0$ if and only if $x=y$;
\item[(b)] $\varrho(y,x)\le C_0\varrho(x,y)$;
\item[(c)] $\varrho(x,y)\le C_1\max\{\varrho(x,z),\varrho(z,y)\}$.
\end{enumerate}
%%%
If $X$ has cardinality at least $2$, then necessarily $C_0,C_1\ge 1$.
A pair $(X,\varrho)$ is called a quasi-metric space. Given $r>0$ and 
$x\in X$, let
\[
B_\varrho(x,r):=\{y\in X\ :\ \varrho(x,y)<r\}
\]
be the quasi-metric ball related to $\varrho$ of radius $r$ and with
center $x$. If $(X,\varrho)$ is a quasi-metric space, then $\mathcal{T}_\varrho$,
the topology on $X$ induced by the quasi-metric $\varrho$ is canonically
defined by declaring $G\subset X$ to be open if and only if for every
$x\in G$, there exists $r>0$ such that $B_\varrho(x,r)\subset G$. The quasi-metric
balls themselves need not be open (unless $\varrho$ is a genuine metric) 
even if $C_0=1$ (see, e.g., an example in \cite[p.~4310]{PS09}).
According to a refined version of the theorem by Mac\'{\i}as and Segovia
(see \cite[Theorem~2]{MS79}) available in \cite[Theorem~2.1]{AM15},
given a quasi-metric $\varrho$, there exists a constant $c\in(0,\infty)$
and a quasi-metric $d$ on $X$ such that for all $x,y\in X$, one has
\[
c^{-1}\varrho(x,y)\le d(x,y)\le c\varrho(x,y),
\quad
d(x,y)=d(y,x),
\]
and all balls $B_d(x,r)$ with respect to $d$ are open in the topology
$\mathcal{T}_\varrho=\mathcal{T}_d$.

From now on, we will assume that $X$ is equipped with this equivalent 
quasi-metric $d$ with the property that all quasi-metric balls $B_d(x,r)$
are open in the topology $\mathcal{T}_d$. For simplicity,
we will write $B(x,r):=B_d(x,r)$.

Let $\mathfrak{M}$ be a $\sigma$-algebra of subsets of $X$ and 
$\mu:\mathfrak{M}\to[0,\infty]$ be a measure. Following 
\cite[Definition~2.9]{AM15}, a measure $\mu$ on the topological space
$(X,\mathcal{T}_d)$ is said to be Borel if $\mathfrak{M}$ contains all
Borel subsets of $X$. A Borel measure $\mu$ on $X$ is said to be
doubling if there exists a constant $A\in(0,\infty)$ such that
\[
0<\mu(B(x,r))\le A\mu(B(x,r/2))<\infty
\]
for all $x\in X$ and $r>0$. In this case the triple $(X,d,\mu)$ is called
the space of homogeneous type.

One says that that a measurable function $f$ on $X$ belongs to
$L_{\rm loc}^1(X,d,\mu)$ if
\[
\int_{B(x,r)}|f(y)|\,d\mu(y)<\infty
\]
for every $x\in X$ and $r>0$.
If $f\in L_{\rm loc}^1(X,d,\mu)$, then the Hardy-Littlewood maximal function
$Mf$ defined by \eqref{eq:Hardy-Littlewood} is measurable on $X$
because the quasi-metric balls $B\subset X$  (with respect to the 
quasi-metric $d$) are open and so $Mf$ is lower semi-continuous.
Further, if $p(\cdot)\in\mathcal{P}(X)$ and
$f\in L^{p(\cdot)}(X,d,\mu)$, then $f\in L_{\rm loc}^1(X,d,\mu)$.
This can be proved as in the Euclidean setting
(see, e.g., \cite[Proposition~2.41]{CF13}).

Following \cite[Definition~3.9]{AM15}, a Borel measure $\mu$ on
$(X,\mathcal{T}_d)$ is said to be Borel-semiregular if for any measurable
set $E$ of finite measure there exists a Borel set $B$ such that 
$\mu(E\Delta B)=0$, where $E\Delta B:=(E\setminus B)\cup(B\setminus E)$.

We will need the following sharp version of the Lebesgue differentiation
theorem (see \cite[Theorem~3.14]{AM15}).
%%%----------------------------------------------------------------------------
\begin{theorem}\label{th:Lebesgue-differentiation-theorem}
Let $(X,d,\mu)$ be a space of homogeneous type. Then the measure $\mu$ is
Borel-semiregular on $(X,\mathcal{T}_d)$ if and only if
\[
\lim_{r\to 0^+}\frac{1}{\mu(B(x,r))}\int_{B(x,r)}f(y)\,d\mu(y)=f(x)
\]
for $\mu$-almost every $x\in X$.
\end{theorem}
%%%----------------------------------------------------------------------------
\section{Proof of the main result}\label{sec:proof} 
Assume that $p_-(X)=1$. Following the general idea of the proof from
\cite[Theorem~3.19]{CF13}, to show that the maximal operator is not bounded, 
we will construct a sequence of functions $\{f_k\}$ such that for all $k$, 
$f_k\in L^{p(\cdot)}(X,d,\mu)$ but the norms $\|Mf_k\|_{p(\cdot)}$ can not 
be uniformly bounded by $\|f_k\|_{p(\cdot)}$.

Since $p_-(X)=1$, for each $k\in\mathbb{N}$ the set 
\[
E_k=\{x\in X\ :\ p(x)<1+1/k\}
\] 
has positive measure. Given that $\mu$ is assumed to be Borel-semiregular, 
applying Theorem~\ref{th:Lebesgue-differentiation-theorem} to the 
function $\chi_{E_k}$, we can choose a point $x_k\in E_k$ such that 
\[
\lim_{r\to0^+}\frac{\mu(E_k\cap B(x_k,r))}{\mu(B(x_k,r))}=1,
\] 
that is, a density point of $E_k$. This choice implies, in particular, that 
for each $k$, there exists a radius $R_k$, $0<R_k<1$, such that if 
$0<r\le R_k$, then 
%%%
\begin{equation}\label{eq:suff-dense-ball}
\frac{\mu(E_k\cap B(x_k,r))}{\mu(B(x_k,r))}>\frac{1+\delta}2,
\end{equation}
%%% 
where $\delta\in(0,1)$ is the reverse doubling constant. 

Let $B_k^0:=B(x_k,R_k)$ be a ball, sufficiently densely---in the sense 
of~\eqref{eq:suff-dense-ball}---filled with the points of $E_k$. For 
$i\in\mathbb{N}_0:=\{0,1,2,\ldots\}$, consider the balls 
\[
B_k^i:=B(x_k,R_k/2^i)
\] 
and split $B_k^0$ into the disjoint union of dyadic annular regions 
$B_k^i\setminus B_k^{i+1}$. Using the doubling property of $\mu$ with 
the constant $A>1$ and the reverse doubling property with the constant 
$\delta$, for each $i\in\mathbb{N}_0$ we estimate beforehand
%%%
\begin{equation}\label{eq:est-meas-annuli}
\mu(B_k^i\setminus B_k^{i+1})\ge(1-\delta)\mu(B_k^i)
\ge
\frac{1-\delta}{A^i}\mu(B_k^0).
\end{equation}
%%%
Finally, define the sequence of functions
%%%
\begin{equation}\label{eq:def-f_k}
f_k(x)=\left(\sum_{i=0}^\infty \frac{\chi_{B_k^i\setminus B_k^{i+1}}(x)}
{A^{i/k}\mu(B_k^i\setminus B_k^{i+1})}\right)\chi_{E_k}(x)
\end{equation}
%%%
on $X$. Note that outside $E_k\cap B_k^0$, the function $f_k$ is identically 
zero.

To show that $f_k\in L^{p(\cdot)}(X,d,\mu)$, we use the simple observation 
that 
\[
f_k(x)^{p(x)}
\le
\max\{1,f_k(x)\}^{p(x)}
\le
\max\{1,f_k(x)^{1+1/k}\}
\le 
1+f_k(x)^{1+1/k}
\]
for all $x\in E_k$, and this, together with~\eqref{eq:est-meas-annuli}, 
gives us 
%%%
\begin{align*}
\rho(f_k)&=\int_{E_k\cap B_k^0} f_k(x)^{p(x)}d\mu(x) \\
&\le \mu(B_k^0) + \int_{B_k^0} f_k(x)^{1+1/k}d\mu(x) \\
&= \mu(B_k^0) + \sum_{i=0}^\infty \frac{\mu((B_k^i\setminus B_k^{i+1})\cap E_k)}
{[A^{i/k}\mu(B_k^i\setminus B_k^{i+1})]^{1+1/k}} \\
&\le \mu(B_k^0) + \sum_{i=0}^\infty \frac{[\mu(B_k^i\setminus B_k^{i+1})]^{-1/k}}{(A^{1/k+1/k^2})^i}\\
&\le \mu(B_k^0) + [(1-\delta)\mu(B_k^0)]^{-1/k} \sum_{i=0}^\infty \frac{A^{i/k}}{(A^{1/k+1/k^2})^i} \\
&= \mu(B_k^0) + [(1-\delta)\mu(B_k^0)]^{-1/k} \sum_{i=0}^\infty \left(A^{-1/k^2}\right)^i,
\end{align*}
%%% 
where the last expression is finite since $A^{-1/k^2}<1$ for each 
$k\in\mathbb{N}$.  

To estimate the norm of $Mf_k$, first fix $x\in E_k\cap B_k^0$. Clearly, 
there exists $i\in\mathbb{N}_0$ such that 
$x\in B_k^i\setminus B_k^{i+1}$ and hence
%%%
\begin{equation}\label{eq:f_k(x)}
f_k(x)=\frac1{A^{i/k}\mu(B_k^i\setminus B_k^{i+1})}.
\end{equation}
%%%
Note that no less than a certain ``portion'' of each annulus 
$B_k^j\setminus B_k^{j+1}$ is filled with the points of 
$E_k$: more precisely, since the radius of each dyadic ball 
$B_k^j$, $j\in\mathbb{N}_0$, does not exceed $R_k$, it 
follows from~\eqref{eq:suff-dense-ball} and the reverse doubling that
%%% 
\begin{align*}
\mu((B_k^j\setminus B_k^{j+1})\cap E_k)&\ge\mu(B_k^j\cap E_k)-\mu(B_k^{j+1}) \\
&> \frac{1+\delta}2\mu(B_k^j)-\delta\mu(B_k^j) \\
&\ge \frac{1-\delta}2 \mu(B_k^j\setminus B_k^{j+1}).
\end{align*}
%%%
Then
%%% 
\begin{align*}
Mf_k(x)&\ge \frac1{\mu(B_k^i)} \int_{B_k^i} f_k(y)\,d\mu(y) 
\\
&= 
\frac1{\mu(B_k^i)} \sum_{j=i}^\infty 
\frac{\mu((B_k^j\setminus B_k^{j+1})\cap E_k)}
{A^{j/k}\mu(B_k^j\setminus B_k^{j+1})} 
\\
&\ge 
\frac1{\mu(B_k^i)} \cdot \frac{1-\delta}2 \sum_{j=i}^\infty (A^{-1/k})^j 
\\
&= 
\frac1{\mu(B_k^i)} \cdot \frac{1-\delta}2 \cdot \frac{A^{-i/k}}{1-A^{-1/k}},
\end{align*}
%%%
which implies, along with~\eqref{eq:est-meas-annuli} and \eqref{eq:f_k(x)}, 
that for $x\in E_k\cap B_k^0$,
\[
Mf_k(x)\ge f_k(x) \cdot \frac{A^{i/k}\mu(B_k^i\setminus B_k^{i+1})}{\mu(B_k^i)} 
\cdot \frac{1-\delta}2 \cdot \frac{A^{-i/k}}{1-A^{-1/k}}
\ge \frac{(1-\delta)^2}{2(1-A^{-1/k})}\,f_k(x).
\]
Trivially, this inequality also holds if $x\not\in E_k\cap B_k^0$. 
Hence, we have shown that 
\[
\|Mf_k\|_{p(\cdot)}\ge\frac{(1-\delta)^2}{2(1-A^{-1/k})} \|f_k\|_{p(\cdot)},
\] 
but since $A^{-1/k}\to1$ as $k\to\infty$, we can not get the uniform 
boundedness of the norms $\|Mf_k\|_{p(\cdot)}$, and this completes the 
proof.
\qed
%%%----------------------------------------------------------------------------
\section{Final remark}\label{sec:final-remark}
If we additionally assume that $p(\cdot)\in\mathcal{P}(X)$ is upper 
semi-continuous, then the hypothesis of Borel-semiregularity of $\mu$
can be dropped because we can avoid using the Lebesgue differentiation
theorem in this case. More precisely, we have the following.
%%%----------------------------------------------------------------------------
\begin{theorem}
Suppose $(X,d,\mu)$ is a space of homogeneous type which has the property that 
the measure $\mu$ is reverse doubling. Given an upper semi-continuous exponent 
function $p(\cdot)\in\mathcal{P}(X)$, if the Hardy-Littlewood maximal operator
$M$ is bounded on the variable Lebesgue space $L^{p(\cdot)}(X,d,\mu)$, then
$p_-(X)>1$.
\end{theorem}
%%%----------------------------------------------------------------------------
\begin{proof}
Assume that $p_-(X)=1$. Since $X$ is open and $p(\cdot)$ is upper 
semi-continuous, for every $k\in\mathbb{N}$, there exist $x_k\in X$ and 
$R_k>0$ such that if $x\in B_k^0:=B(x_k,R_k)$, then $p(x)<1+1/k$.
Now define $f_k$ replacing $\chi_{E_k}$ by $\chi_{B_k^0}$ in \eqref{eq:def-f_k}.
After this the proof goes as that of Theorem~\ref{th:main} with minor changes. 
\end{proof}
%%%----------------------------------------------------------------------------
\section*{Acknowledgments}
We would like to thank David Cruz-Uribe for useful discussions and for sharing 
with us the thesis of Roberts \cite{R23} written under his guidance.
%%%----------------------------------------------------------------------------
\bibliographystyle{abbrv}
\bibliography{SHT-bib}

\begin{thebibliography}{10}

\bibitem{AM15}
R.~Alvarado and M.~Mitrea.
\newblock {\em Hardy spaces on {A}hlfors-regular quasi metric spaces}, volume
  2142 of {\em Lecture Notes in Mathematics}.
\newblock Springer, Cham, 2015.
\newblock A sharp theory.

\bibitem{CFN03}
D.~Cruz-Uribe, A.~Fiorenza, and C.~J. Neugebauer.
\newblock The maximal function on variable {$L^p$} spaces.
\newblock {\em Ann. Acad. Sci. Fenn. Math.}, 28(1):223--238, 2003.

\bibitem{CF13}
D.~V. Cruz-Uribe and A.~Fiorenza.
\newblock {\em Variable {L}ebesgue spaces}.
\newblock Applied and Numerical Harmonic Analysis. Birkh\"{a}user/Springer,
  Heidelberg, 2013.
\newblock Foundations and harmonic analysis.

\bibitem{DFGL08}
G.~Di~Fazio, C.~E. Guti\'{e}rrez, and E.~Lanconelli.
\newblock Covering theorems, inequalities on metric spaces and applications to
  {PDE}'s.
\newblock {\em Math. Ann.}, 341(2):255--291, 2008.

\bibitem{DHHMS09}
L.~Diening, P.~Harjulehto, P.~H\"{a}st\"{o}, Y.~Mizuta, and T.~Shimomura.
\newblock Maximal functions in variable exponent spaces: limiting cases of the
  exponent.
\newblock {\em Ann. Acad. Sci. Fenn. Math.}, 34(2):503--522, 2009.

\bibitem{DHHR11}
L.~Diening, P.~Harjulehto, P.~H\"{a}st\"{o}, and M.~R\r{u}\v{z}i\v{c}ka.
\newblock {\em Lebesgue and {S}obolev spaces with variable exponents}, volume
  2017 of {\em Lecture Notes in Mathematics}.
\newblock Springer, Heidelberg, 2011.

\bibitem{INS13}
M.~Izuki, E.~Nakai, and Y.~Sawano.
\newblock The {H}ardy-{L}ittlewood maximal operator on {L}ebesgue spaces with
  variable exponent.
\newblock In {\em Harmonic analysis and nonlinear partial differential
  equations}, RIMS K\^{o}ky\^{u}roku Bessatsu, B42, pages 51--94. Res. Inst.
  Math. Sci. (RIMS), Kyoto, 2013.

\bibitem{ISN14}
M.~Izuki, E.~Nakai, and Y.~Sawano.
\newblock Function spaces with variable exponents---an introduction---.
\newblock {\em Sci. Math. Jpn.}, 77(2):187--315, 2014.

\bibitem{MS79}
R.~A. Mac\'{\i}as and C.~Segovia.
\newblock Lipschitz functions on spaces of homogeneous type.
\newblock {\em Adv. in Math.}, 33(3):257--270, 1979.

\bibitem{M11}
L.~Maligranda.
\newblock Hidegoro {N}akano (1909--1974)---on the centenary of his birth.
\newblock In {\em Banach and function spaces {III} ({ISBFS} 2009)}, pages
  99--171. Yokohama Publ., Yokohama, 2011.

\bibitem{PS09}
M.~Paluszy\'{n}ski and K.~Stempak.
\newblock On quasi-metric and metric spaces.
\newblock {\em Proc. Amer. Math. Soc.}, 137(12):4307--4312, 2009.

\bibitem{R23}
T.~Roberts.
\newblock Necessary conditions for bounded fractional maximal operators and
  fractional singular integral operators on variable {L}ebesgue spaces.
\newblock In {\em Master of Arts Thesis}. The University of Alabama,
  Tuscaloosa, Alabama, 2023.

\end{thebibliography}
\end{document}